\documentclass[11pt]{amsart}
\usepackage{geometry}               
\usepackage{graphicx}
\usepackage{amssymb}
\usepackage[colorlinks=true,urlcolor=blue]{hyperref}

\newtheorem{theorem}{Theorem}[section]
\newtheorem{lemma}[theorem]{Lemma}

\theoremstyle{definition}

\newtheorem{proposition}[theorem]{Proposition}
\newtheorem{corollary}[theorem]{Corollary}

\usepackage{lineno}

\title{Conditional probability of derangements and fixed points}
\author{Sam Gutmann, Mark Mixer, and Steven Morrow}

\begin{document}

\maketitle

\begin{abstract}
The probability that a random permutation in $S_n$ is a derangement is well known to be $\displaystyle\sum\limits_{j=0}^n (-1)^j \frac{1}{j!}$.  In this paper, 
we consider the conditional probability that the $(k+1)^{st}$ point is fixed, given there are no fixed points in the first $k$ points.
We prove that when $n \neq 3$ and $k \neq 1$, this probability is a decreasing function of both $k$ and $n$.   Furthermore, it is proved that this conditional probability is well approximated by $\frac{1}{n} - \frac{k}{n^2(n-1)}$.    Similar results are also obtained about the more general conditional probability that the $(k+1)^{st}$ point is fixed, given that there are exactly $d$ fixed points in the first $k$ points.



\end{abstract}

\section{Introduction}
In the famous hat-check problem, the reader is asked to suppose there are $n$ people at a party, each of whom brought a hat; at the end of the party everyone picks up a random hat.  It is well known that the probability that no one receives their own hat is the same as the probability that a random permutation has no fixed points.
 The classical solution (see for example~\cite{Feller}) is to write:
\begin{equation}\label{eq1}
P(A_1^c \cap A_2^c \cap \dots \cap A_n^c)  = 1 - P(A_1 \cup \dots \cup A_n) = \displaystyle\sum\limits_{j=0}^n (-1)^j \frac{1}{j!} \approx 1/e.
\end{equation}
The solution can be denoted by $d_n/n!$, where $d_n$ is the number of derangements (i.e. permutations with no fixed points) of a set of sized $n$.  In~\cite{Takacs}, there is a compendium of the history of the topic of fixed points of permutations.   Additionally, many variations of fixed point questions have been studied (see for instance: \cite{Fisk, HSW,  Penrice, Rawlings}).  

In the present paper, we consider the conditional probability of continuing, or not continuing, a derangement in-progress; i.e., if we know there are no fixed points among the first $k$ points (a ``partial derangement''), we look at the chance that the next point is fixed, defining a function $f$ for this conditional probability. We prove that this function decreases as $k$ increases, except when $n=3$.  
In the hat-check problem, this means that with more than 3 people, the chance a hat is returned to its owner gets worse as more hats are returned unsuccessfully.

We also generalize by dealing with permutations with some fixed points, where we examine the behavior of the function $f(n,k,d)$ which gives the probability that $(k+1)$ is fixed, given there are exactly $d$ fixed points in the first $k$ points.  This is aided by a recursion relating this probability to a partial derangement.   From this recursion and the inclusion-exclusion principle, one can show that 
\begin{equation} \label{eqfnkd1}
f(n,k,d) =  1 -  \frac{ \displaystyle \sum_{j=0}^{k-d+1} (-1)^j \binom{k-d+1}{j} (n-d-j)!}{ \displaystyle \sum_{j=0}^{k-d} (-1)^j \binom{k-d}{j} (n-d-j)!}.
\end{equation}
or equivalently that 
\begin{equation} \label{eqfnkd2}
f(n,k,d) =  \frac{ \displaystyle \sum_{j=0}^{k-d} (-1)^j \binom{k-d}{j} (n-d-j-1)!}{ \displaystyle \sum_{j=0}^{k-d} (-1)^j \binom{k-d}{j} (n-d-j)!}.
\end{equation}

However, neither of these formulas provides much insight into the behavior of $f$.    We will produce an approximation for this conditional probability $f(n,k,d)$, and use it to determine the monotone behavior of $f$ (with a few noted exceptions) in each of its variables.  
There is related work, such as~\cite{FDS},
which shows that the set of fixed points of a random permutation has the `FKG' property, meaning that any two reasonable definitions of what it means for the set of fixed points to be `big' are positively correlated. The present paper looks at a similar problem, but our results do not follow from the FKG property.

\begin{figure}[h]
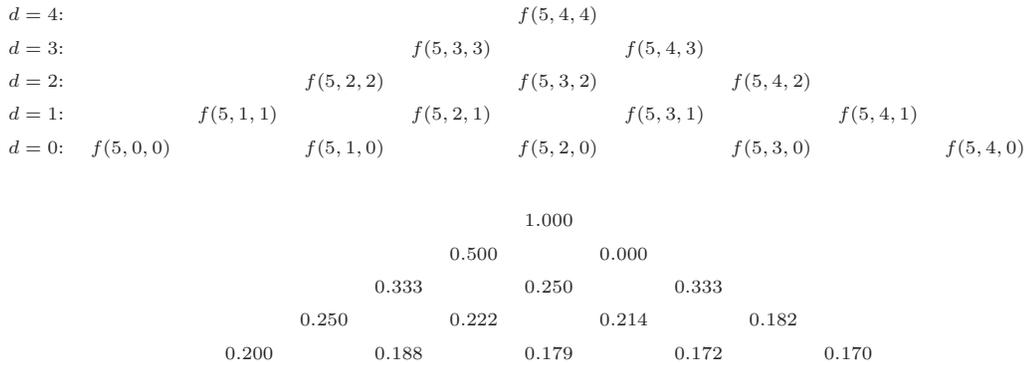

\tiny
\begin{center}
\begin{tabular}{rccccccccc}
$d=4$:&    &    &    &    &  $f(5,4,4)$\\\noalign{\smallskip\smallskip}
$d=3$:&    &    &    &  $f(5,3,3)$ &    &  $f(5,4,3)$\\\noalign{\smallskip\smallskip}
$d=2$:&    &    &  $f(5,2,2)$ &    &  $f(5,3,2)$ &    &  $f(5,4,2)$\\\noalign{\smallskip\smallskip}
$d=1$:&    &  $f(5,1,1)$ &    &  $f(5,2,1)$ &    &  $f(5,3,1)$ &    &  $f(5,4,1)$\\\noalign{\smallskip\smallskip}
$d=0$:&  $f(5,0,0)$ &    &  $f(5,1,0)$ &    &  $f(5,2,0)$ &    &  $f(5,3,0)$ &    &  $f(5,4,0)$\\\noalign{\smallskip\smallskip}
\end{tabular}

\vspace{.5cm}

\begin{tabular}{rccccccccc}
\hspace{.5cm} &    &    &    &    &  $1.000$\\\noalign{\smallskip\smallskip}
 &   &    &    &  $0.500$ &    &  $0.000$\\\noalign{\smallskip\smallskip}
  &  &    &  $0.333$ &    &  $0.250$ &    &  $0.333$\\\noalign{\smallskip\smallskip}
   & &  $0.250$ &    &  $0.222$ &    &  $0.214$ &    &  $0.182$\\\noalign{\smallskip\smallskip}
 & $0.200$ &    &  $0.188$ &    &  $0.179$ &    &  $0.172$ &    &  $0.170$\\\noalign{\smallskip\smallskip}
\end{tabular}
\caption{A triangle $f(5,k,d)$, and its values rounded to 3 decimal places \label{fign=5} }
\end{center}

\end{figure}
\normalsize

To help illustrate our results, in Figure~\ref{fign=5}, we construct a triangle as a visualization for these conditional probabilities $f(n,k,d)$ when $n=5$.  We will show that except in one case, $f(n,k,d)$ decreases in $k$. In the triangle, each horizontal row from left to right decreases, except for the row of three entries that occurs when $d=n-3$.  Also, we will show that except in one case $f(n,k,d)$ decreases in $n$.   After applying a recursion, this will mean that in the triangle each `tilted' column running from top right to bottom left decreases, except the column which begins with $f(n,n-1,n-2)= 0$.  Similarly, we will show that except in one case $f(n,k,d)$ increases in $d$.  In the triangle, each `tilted' column running from top left to bottom right decreases, again excepting the column that contains $f(n,n-1,n-2)= 0$.  

The remaining sections of the paper are organized as follows.   In Section~\ref{secBasics} one will find notation, background, and some basic results.
In Section~\ref{secBijections}, in order to formalize some counting arguments in our proofs, we analyze various subsets of $S_n$ based on their fixed points.
Then in Section~\ref{secRecursive}, we provide recursive results that allow the reduction to partial derangements.   Most of the technical work is done in Section~\ref{dZero},  where we prove some results about partial derangements. Finally, Section~\ref{secResults} contains the proof of the theorems, which give an approximation for $f(n,k,d)$, and show its monotone behavior in each of its variables.  This is accomplished by applying the recursive results of  Section~\ref{secRecursive} to the main results of Section~\ref{dZero}.


\section{Background and Statement of Theorems}
\label{secBasics}

\subsection{Notation}
Let $[n]$ denote the collection of integers $[1,\ldots, n]$.  
$S_n$ is the symmetric group acting on $[n]$.  It will be useful to think of $S_n$ simply as the set of bijections from $[n]$ to itself: $S_n= \{ \alpha \ | \ \alpha : [n] \rightarrow [n]$\}.
For any permutation $\alpha \in S_n$ we will be interested in its fixed points.  Let $F(\alpha) = \{ i : \alpha(i) = i \}$.  
Define $S_{n,k,d}$ as the collection of permutations in $S_n$ with exactly $d$ fixed points in the first $k$ points.  That is to say $$S_{n,k,d}:= \{ \alpha \in S_n \ : \  |F(\alpha) \cap [k] | = d\}.$$     We can use $S_{n,k,d}$ to define two functions of interest.  Let $c(n,k,d) := | S_{n,k,d} |$ and $p(n,k,d)=c(n,k,d)/n!$ which is the probability that a random element of $S_n$ has exactly $d$ fixed points in the first $k$ points.  

\subsection{Basic Notions and Statement of Results}

When $k=n$ and $d=0$, $c(n,n,0)$ equals the number of derangements of a set of size $n$, often denoted $d_n$; in this case, $p(n,n,0)$ is given in Equation~\ref{eq1}.  
For other values of $k$, using the inclusion-exclusion principle, it is straightforward to show:

\begin{equation} \label{eqcnk0}
c(n,k,0) = \sum_{j=0}^{k} (-1)^j \binom{k}{j} (n-j)!.
\end{equation}

This is a count of $k$-partial derangements, those permutations with no fixed points among the first $k$ points. 
It should be noted that our use of the term ``partial derangement" differs from the usage in~\cite{Antonelli} and~\cite{EHH} where the expression is used to describe any permutation which has at least one fixed point and at least one non-fixed point. 

Now assume that it is known that there are exactly $d$ fixed points in the first $k$ points, we are interested in the conditional probability that $(k+1)$ is fixed.  Call this conditional probability $f(n,k,d)$.  
In fact, the values of $p(n,k,d)$ and $f(n,k,d)$ do not depend on the arbitrary choice to examine the first $k$ points.  We can define the functions more generally as follows.

Let $A$ be any subset of $[n]$ where $|A| = k$, and let $a$ be any element of $[n] \setminus A$.  
Then the value of $p(n,k,d) := P(\text{exactly }d \text{ points in } A \text{ are fixed}) $ is independent of your choice for $A$, since the points can be relabeled to be the first $k$ points accordingly.    Thus
\begin{equation} \label{eqp0}
p(n,k,d) = P(\{ \alpha \in S_n \ : \  |F(\alpha) \cap A | = d\}).
\end{equation}

Similarly, $f(n,k,d)$ is the conditional probability that $a$ is fixed, assuming there are exactly $d$ fixed points in $A$.  This is also independent of your choices for $A$ and $a$, and thus
\begin{equation} \label{eqf0}
 f(n,k,d) = P( \{ \alpha \in S_n \ : \    \alpha(a)=a \} | \{  \alpha \in S_n \ : \  |F(\alpha) \cap A | = d \} ).
\end{equation}

When the set $A$ and the point $a$ are not given, we will assume them to be $[k]$ and $(k+1)$ respectively.  We note here that we allow $k=0$, where $p(n,0,0)=1$, $c(n,0,0)=n!$, and $f(n,0,0) = P( \{ \alpha \in S_n \ : \    \alpha(a)=a \}) = \frac{1}{n}$.
The parameters $n,k,d$ are all integers where $n\geq 1$, and $0 \leq d \leq k$.  For $p(n,k,d)$ we allow $0\leq k \leq n$, and for $f(n,k,d)$ we require that $0\leq k \leq n-1$.

We now provide a straightforward lemma and then derive an expression for $f(n,k,d)$.

\begin{lemma} \label{lemfixed1}
Let $x$ and $k$ be in $[n]$.  If $x>k$ then, the number of permutations of $[n]$ that both have $d$ fixed points in the first $k$ points and have $x$ fixed is equal to the number of permutations of $[n-1]$ that have $d$ fixed points in the first $k$ points.     If $x \leq k$, the number of permutations of $[n]$ that both have $d$ fixed points in the first $k$ points and have $x$ fixed is equal to the number of permutations of $[n-1]$ that have $d-1$ fixed points in the first $k-1$ points.
\end{lemma}

This can be proved by removing the point $x$ from $[n]$ and relabeling the points accordingly.   In our notation:

$$ | \{ \alpha \in S_n \ : \  |F(\alpha) \cap [k] | = d , \  \alpha (x) = x  \} | = c(n-1,k,d) \text{ when } x > k,$$
$$ | \{ \alpha \in S_n \ : \  |F(\alpha) \cap [k] | = d , \  \alpha (x) = x  \} | = c(n-1,k-1,d-1) \text{ when } x \leq k.$$

In particular, when considering the function $f(n,k,d)$, it will be useful for us to notice that, for a random permutation in $S_n$, $$P(d \text{ of the first } k \text{ points and } (k+1) \text{ are fixed}) = \frac{ c(n-1,k,d)}{n!}.$$

 The need for a more general version of this result, where the image of a point is known but it need not be fixed, leads to the work in Section~\ref{secBijections}.

From Lemma~\ref{lemfixed1} it follows that, for $0 \leq k < n,$ $0 \leq d \leq k$:

\begin{equation} \label{eqfnkd}
f(n,k,d)= \frac{c(n-1,k,d)}{c(n,k,d)}.  
\end{equation}

\subsection{Special Cases}

Much of our work will focus on the case when $d=0$.  
There exist some well-known recursive results on the number of derangements $d_n = c(n,n,0)$. Let us define $d_{n,k} := c(n,k,0)$ to be the number of $k$-partial derangements on $[n]$. Then: 
\begin{equation} \label{eqsub}
d_{n,n-1}=d_n + d_{n-1} \mbox{ for } n\geq 1. 
\end{equation}
To prove this, consider the set $S_{n,n-1,0}$, whose cardinality equals $d_{n,n-1}$.  This set contains all the permutations that have no fixed points in the first $(n-1)$ points and is the disjoint union of two subsets: those with and those without the last point fixed. The subset without the final point fixed has cardinality $d_n$, and the subset with the final point fixed has cardinality $d_{n-1}$.  A similar argument can be used to prove a more general version of this recurrence: 
\begin{equation} \label{eqsub2}
d_{n,k}=d_{n,k+1} + d_{n-1,k} \mbox{ for } 0 \leq k \leq n-1.
\end{equation}

We are able to create new expressions for the number of derangements $d_n$.  For example, using straightforward combinatorial arguments, we can relate $d_n$ to numbers of permutations with exactly 1 fixed point. We have that $(n+1) d_n =c(n+1,n+1,1)$. Perhaps less obvious is the relation:
\begin{equation} \label{eqdncn}
d_n = c(n,n-1,1).
\end{equation}
This is true by way of making a bijection between two subsets of $S_n$: $S_{n,n,0}$ (the set of derangements on $[n]$) and $S_{n,n-1,1}$.  For each derangement $\sigma$, we know $\sigma(n) = i \neq n$. Map $\sigma$ to the permutation $\tau$ which is identical on all points except with $\tau(i)=i$ and $\tau(n)=\sigma(i)$. That is, we interchange the images of point $n$ and point $\sigma(n)$. $\tau$ is then guaranteed to have exactly one fixed point among the first $n-1$ points, with either a fixed or non-fixed point in the last point, $n$.  The set of all such permutations makes up precisely the set $S_{n,n-1,1}$.

We are also able to provide a new proof of the well-known recurrence $d_n = (n-1)(d_{n-1} + d_{n-2})$, as follows. Again consider the set $S_{n,n-1,1}$, and partition it into subsets according to which of the first $n-1$ points is fixed.  Each of these subsets has cardinality $c(n-1,n-2,0)$, by 
Lemma~\ref{lemfixed1}.
Since there are $n-1$ such subsets, it follows that $|S_{n,n-1,1}| = c(n,n-1,1) = (n-1)c(n-1,n-2,0)$. Combining this with Equations~\ref{eqsub} and \ref{eqdncn} produces the desired result.

When $d=0$, by making use of Equation~\ref{eqsub2}, we also have:
\begin{equation} \label{eqfnk0b}
f(n,k,0) = 1 - \frac{c(n,k+1,0)}{c(n,k,0)} = 1 - \frac{p(n,k+1,0)}{p(n,k,0)}.
\end{equation}

Using this with Equation~\ref{eqcnk0} gives a way to compute values of $f(n,k,0)$.

\begin{equation} \label{eqftop}
f(n,k,0) = 1 - \frac{p(n,k+1,0)}{p(n,k,0)}  = 1 - \frac{\displaystyle \sum_{j=0}^{k+1} (-1)^j \binom{k+1}{j} (n-j)!}{\displaystyle \sum_{j=0}^{k} (-1)^j \binom{k}{j} (n-j)!}
\end{equation}

When $k$ is small we can algebraically simplify the function $f(n,k,0)$.  
For instance, when $k=0$, $f(n,0,0)$ is the probability that the first point of a random permutation in $S_n$ is fixed. Thus $f(n,0,0) = \frac{1}{n}$. 
 Table~\ref{tab3} gives simplifications of the function $f(n,k,0)$ for $k \leq 3$.

\begin{table}[h]
\begin{center}
\begin{tabular}{|c|c|c|c|c|}
\hline
$k$ & 0 & 1 & 2 & 3 \\ \hline 
$f(n,k,0)$  &  $\frac{1}{n}$ & $ \frac{(n-2)}{(n-1)^2}$ & $  \frac{n^2-5n+7}{(n-2)(n^2-3n+3)}$ & $\frac{n^3-9n^2+29n-34}{(n-3)(n^3-6n^2+14n-13)}$ \\ \hline
\end{tabular}
\caption{Simplifications of $f(n,k,0)$ for small values of $k$.  \label{tab3} }
\end{center}
\end{table}

Some small values of $p(n,k,0)$ and $f(n,k,0)$ will be useful in later proofs.
Tables~\ref{tab1}~and~\ref{tab2} provide some values of  $p(n,k,d)$ and $f(n,k,d)$ for $d=0$ and small $k$ and $n$.  The case where $d >0$ will be handled in Section~\ref{secRecursive}.

\begin{table}[h]
\begin{center}

\begin{tabular}{|c|c|c|c|c|c|c|c|}
\hline
$n \backslash k$ & 0 & 1 & 2 & 3 & 4 & 5 & 6 \\ \hline
1 & 1  & 0 &&&&& \\ \hline
2 & 1 & .5 & .5 &&&&\\ \hline 
3 & 1 & .6666 & .5 & .3333 &&& \\ \hline
4 & 1 & .75 & .5833 & .4583 & .375 && \\ \hline
5 & 1 & .8 & .65 & .5333 & .4417 & .3667 & \\ \hline
6 & 1 & .8333 & .7 & .5917 & .5203 & .4292 & .3681 \\\hline
\end{tabular}
\caption{$p(n,k,0)$ for small values of $k,n$ (rounded to 4 decimal places) \label{tab1}}
\end{center}

\end{table}

\begin{table}[h]
\begin{center}

\begin{tabular}{|c|c|c|c|c|c|c|c|}
\hline
$n \backslash k$ & 0 & 1 & 2 & 3 & 4 & 5 \\ \hline
1 & 1 &&&&&\\ \hline
2 & .5 & 0 &&&&\\ \hline
3 & .3333 & .25 & .3333 &&&\\ \hline
4 & .25 & .2222 & .2143 & .1818 &&\\ \hline
5 & .2 & .1875 & .1795 & .1719 & .1698 & \\ \hline
6 & .1667 & .16 & .1548 & .1502 & .1464 & .1424\\ \hline
\end{tabular}
\caption{$f(n,k,0)$ for small values of $k,n$ (rounded to 4 decimal places) \label{tab2}}
\end{center}

\end{table}

From these tables, one notices some patterns. 
For instance, $p(n,k,0)$ looks to be nonincreasing in $k$; this is shown in Proposition~\ref{prop1}.
Also $p(n,k,0)$ looks to be nondecreasing in $n$; this is shown in Proposition~\ref{prop3}.  
The values in Table~\ref{tab2} are the same as the values in Figure~\ref{fign=5}.  This leads to a recursion shown in Proposition~\ref{proprecursion}.

Those results are generally straightforward; however some results about $f(n,k,0)$ are more technical.
For instance, for $n \neq 3$, $f(n,k,0)$ looks to be decreasing in $k$; this is shown in Theorem~\ref{thm2}.
Finally, for $k \neq 1$, $f(n,k,0)$ looks to be decreasing in $n$; this is shown in Theorem~\ref{thmnewn}.
These results will be utilized to show the monotone behavior of $f(n,k,d)$ in each of its variables (with noted exceptions) in Section~\ref{secResults}.


\section{Bijections}
\label{secBijections}

In this section we will look at subsets of $S_{n,k,d}$, where the image of one point $i$ is known.   These subsets will provide a useful tool to aid in understanding the function $f(n,k,d)$.
Let $B_{n,i,j}=   \{ \rho \ | \ \rho \in S_n, \ \rho(i)=j \}  $ and $B_{n,k,d,i,j} = S_{n,k,d} \cap B_{n,i,j}$.  The structure of $B_{n,k,d,i,j}$ will depend on the values of $i$ and $j$. 

To examine each $B_{n,k,d,i,j}$ closely, we provide an explicit bijection $\Psi$ between $B_{n,i,j}$ and $S_{n-1}$ which allows precise tracking of fixed points.   The bijection $\Psi$ will be constructed as the composition of either two or three other bijections  (depending on if $i = j$). Consider the following other sets of bijective functions, with differing domains and ranges depending on $i$ and $j$.

Let $B'_{n,i,j} :=   \{ \rho \ | \ \rho: ([ n ] \setminus i)  \rightarrow ([ n ]  \setminus j)  \}  $.  Note that it is possible for $i=j$ and where we consider the following set of bijections $B'_{n,i,i} :=   \{ \rho \ | \ \rho: ([ n ] \setminus i)  \rightarrow ([ n ]  \setminus i)  \}  $.

We can construct a bijection from $B_{n,i,j}$ to $S_{n-1}$ by showing that the following maps $\Psi', \Psi'',$ and $\Psi'''$ are bijections.

$\Psi' : B_{n,i,j} \rightarrow B'_{n,i,j}$.   $\Psi'(\rho)$ is the restriction of $ \rho$ to $ [ n ] \setminus i$. So $\Psi'(\rho) = \alpha$ where $\alpha(x) = \rho(x)$ for $x \neq i$.

$\Psi'' : B'_{n,i,j} \rightarrow B'_{n,i,i}$.   $\Psi''(\rho) = \alpha$ where: 

$$
\alpha(x) =\begin{cases}
\rho(x), & \text{if }\rho(x) \neq i\\
j , & \text{if } \rho (x) = i	 \end{cases}
$$

$\Psi''' : B'_{n,i,i} \rightarrow S_{n-1}$.   $\Psi'''$ is defined by relabeling the points in the domain and range of $\rho$.  If a point is larger than $i$, we subtract 1.  Formally,

$\Psi'''(\rho) = \alpha$ where:

$$
\alpha(x) = \begin{cases}
\rho(x) & \text{if } x, \rho(x) < i \\ 
\rho(x) - 1& \text{if } x < i, \rho(x) > i \\ 
\rho(x+1 )& \text{if } x \geq i,  \rho(x) < i \\
\rho(x+1 ) - 1& \text{if } x, \rho(x) \geq  i  
\end{cases}
$$

 Composing these maps gives us $\Psi: B_{i,j} \rightarrow S_{n-1}$. Note we only consider $\Psi''$ when $i \ne j$, as when $i=j$, $\Psi''$ simplifies to the identity mapping.   
Furthermore, keeping track of the fixed points of these maps will be essential.
To illustrate these mappings, we provide the following example where $i=4$ and $j=5$.  Let $\rho \in S_8$.  

$$\rho = \begin{bmatrix} 1 & 2 & 3 & 4 & 5 & 6 & 7 & 8 \\ 7 & 2 & 1 & 5 & 8 & 6 & 4 & 3 \end{bmatrix}.$$

$$\Psi'(\rho) = \begin{bmatrix} 1 & 2 & 3  & 5 & 6 & 7 & 8 \\ 7 & 2 & 1 & 8 & 6 & 4 & 3 \end{bmatrix}.$$

$$\Psi''(\Psi'(\rho)) = \begin{bmatrix} 1 & 2 & 3  & 5 & 6 & 7 & 8 \\ 7 & 2 & 1 & 8 & 6 & 5 & 3 \end{bmatrix}.$$
 
 Finally, $$\Psi(\rho):=\Psi'''(\Psi''(\Psi'(\rho))) = \begin{bmatrix} 1 & 2 & 3  & 4 & 5 & 6 & 7 \\ 6 & 2 & 1 & 7 & 5 & 4 & 3 \end{bmatrix}.$$

\begin{lemma} \label{lem0}
$\Psi'$, $\Psi''$, and $\Psi'''$ are bijections.  Furthermore, the following is true about fixed points.  
\begin{enumerate}
\item Let $\rho \in B_{n,i,j}$.  $\rho(a)=a$ if and only if $(\Psi' \rho)(a)=a$ or $i=j$. 
\item Let $i\neq j$ and $\rho \in B'_{n,i,j}$. If $(\Psi'' \rho) (a) = a$, then either $\rho(a)=a$ or $\rho(a)=i$.
\item Let $\rho \in B'_{n,i,i}$. When $a < i$ then $\rho(a)=a$ if and only if $(\Psi''' \rho) (a)=a$.  When  $a > i$ then $\rho(a)=a$ if and only if $(\Psi''' \rho) (a-1)=a-1$
\end{enumerate}

\end{lemma}
\begin{proof}
This fact that $\Psi'$, $\Psi''$, and $\Psi'''$ are bijections follows directly from their definitions as each map is defined to be surjective and easily reversible. Let $\rho \in B_{n,i,j}$ and assume $\rho(a)=a$.  If $a \ne i$ then $(\Psi' \rho)(a)=a$.     Conversely, since $j=\rho(i)$, if $(\Psi' \rho)(a)=a$ or $i=j$, then $\rho(a)=a$.
Next let $\rho \in B'_{n,i,j}$ where $i \neq j$ and assume $(\Psi'' \rho) (a) = a$.   If $a=j$, then by definition $\rho(a)=i$.    If $a \neq j$ then $\rho(a) = a$.
Finally, let $\rho \in B'_{n,i,i}$.  Assume $\rho(a)=a$ where $a< i$. Then since $a, \rho(a) < i$, $(\Psi''' \rho)(a) = \rho(a) =a.$ On the other hand, if $a > i$ then $(\Psi''' \rho) (a-1)= \rho(a-1+1) -1 = a -1$.
\end{proof}

The bijections in Lemma~\ref{lem0} allow us to determine the size of useful subsets of $B_{i,j}$ based on fixed points.
 In each of the corollaries below we assume that $n \geq 2$, $1 \leq k \leq n-1$, $1 \leq d \leq k-1$, and $i,j \in [n]$.

\begin{corollary}\label{cor_kij}
   If $k  < i,j $, then $| B_{n,k,d,i,j} | =  c(n-1,k,d)$.
\end{corollary}
\begin{proof}
Let $\rho \in B_{n,k,d,i,j}$ where $k  < i,j $.  Then $\Psi \rho$ has exactly $d$ fixed points in the first $k$ points, regardless of whether $i=j$.
\end{proof}
We note here that when $i=j$, this is also proved in Lemma~\ref{lemfixed1}.

\begin{corollary}\label{cor_ikj}
   If $ i \leq k < j $, then $| B_{n,k,d,i,j} | =  c(n-1,k-1, d)$.
\end{corollary}
\begin{proof}
Let $\rho \in B_{n,k,d,i,j}$ where $ i \leq k < j $.  Then $\Psi \rho$ has exactly $d$ fixed points in the first $k-1$ points.
\end{proof}
In this last case, we note that when $k=1$, then $d=0$ as $i \leq k = 1$ is not fixed.  Thus $c(n-1,k-1, d) = |S_{n-1,0, 0}|=(n-1)!$.


\section{Recursive Results}
\label{secRecursive}

We start this section with a technical lemma about $p(n,k,d)$.

\begin{lemma} \label{lem2}
For $0 < k < n$, $$p(n,k,d) = \frac{k}{n} p(n-1,k-1,d) + (1-\frac{k}{n})p(n-1,k,d).$$
\end{lemma}
\begin{proof}
Recall that $p(n,k,d) = c(n,k,d)/n!$.  
Let $\alpha \in S_{n,k,d}$.  We break $p(n,k,d)$ into two options based on what point gets sent to $n$.  
For the $k$ cases when $\alpha(i) = n$ for $i\leq k$, we know that $\alpha \in B_{n,k,d,i,n}$, and by Corollary~\ref{cor_ikj}, there are $c(n-1,k-1,d)$ such permutations.  Thus the probability of this option is $k\frac{c(n-1,k-1,d)}{n!}$ .
On the other hand, for the $n-k$ cases when $\alpha(i) = n$ for $i> k$, by Corollary~\ref{cor_kij}, there are $c(n-1,k,d)$ such permutations.  Thus, the probability of this option is $(n-k)\frac{c(n-1,k,d)}{n!}$.  Rewriting we get 
\begin{align*}
p(n,k,d) = \frac{k}{n} \frac{c(n-1,k-1,d)}{(n-1)!} + \frac{(n-k)}{n}\frac{c(n-1,k,d)}{(n-1)!} \\
= \frac{k}{n}  p(n-1,k-1,d) + (1-\frac{k}{n})p(n-1,k,d).
 \end{align*}
\end{proof}

To count $c(n,k,d)$ when $d$ is not 0, for each particular choice of which $d$ points are fixed, there are $k-d$ non-fixed points to be filled by $n-d$ choices, and therefore $c(n-d,k-d,0)$ ways to do so.  Therefore, we have:

\begin{equation}\label{countd}
c(n,k,d) = \binom{k}{d} c(n-d,k-d,0).
\end{equation}

Combining Equations~\ref{eqcnk0} and~\ref{countd}, we get the following expression for $c(n,k,d)$, which also yields an expression for $p(n,k,d) = \frac{c(n,k,d)}{n!}$.

\begin{equation}\label{formc}
c(n,k,d) =  \displaystyle \sum_{j=0}^{k-d} (-1)^j \binom{k-d}{j} (n-d-j)! \binom{k}{d}
\end{equation}

Using Equation~\ref{countd}, it follows that 

\begin{equation}\label{recursion}
c(n,k,i) = \frac{\binom{k}{i}}{\binom{k-j}{i-j}} c(n-j,k-j,i-j), ~\mbox{ ~for } 0 \leq j \leq i.
\end{equation}

We can then use this recursion for $c(n,k,d)$ to provide a recursion for $f(n,k,d)$.

\begin{proposition} \label{proprecursion}
$f(n,k,d) = f(n-a, k-a, d-a)$ for $0 \leq k < n, 0 \leq d \leq k$, and for all $a$ where $0 \leq a \leq d$.
\end{proposition}
\begin{proof}
\begin{align*}
f(n,k,d) & = \frac{c(n-1,k,d)}{c(n,k,d)} \text{ ~(by Equation~\ref{eqfnkd})} \\
& = \frac{c(n-1-a,k-a,d-a)}{c(n-a,k-a,d-a)} \text{ ~(by Equation~\ref{recursion}, after canceling the binomial coefficients)} \\
& = f(n-a,k-a,d-a) \text{ ~(by Equation~\ref{eqfnkd})}.
\end{align*}
\end{proof}

To conclude this section, we demonstrate how our results will utilize the fact that $f(n,k,d) = f(n-d, k-d, 0)$.    For example, in Figure~\ref{fign=5}, only the last horizontal row has $d=0$.  Using Proposition~\ref{proprecursion}, this same triangle of values can thus be rewritten as a triangle with different inputs for $f$ as shown in Figure~\ref{fign=6}. In other words, for each $n$, we can equate such a triangle like Figure~\ref{fign=5} to one where each row consists of values of $f$ for partial derangements for smaller values of $n$. Therefore, the values for all triangles constructed as in Figure~\ref{fign=5} for some $n$ are identical except for having different numbers of rows. 

Notice the first entry in each row is always $1/n$.  The values in each row then decrease (except when $n=3$), ending with a value very nearly equal to $1/(n+1)$, by way of the estimate $f(n,n-1,0) \approx \frac{1}{n+1}$ for large $n$. This is derived by combining Equation~\ref{eqsub} with the well-known result $d_n = n d_{n-1} + (-1)^n$, and then making use of the relation $f(n,k,0)=\frac{c(n-1,k,0)}{c(n,k,0)} = \frac{d_{n-1,k}}{d_{n,k}}$ from Equation~\ref{eqfnkd}.  This yields $f(n,n-1,0) = \frac{d_{n-1}}{d_{n,n-1}} = \frac{d_{n-1}}{(n+1)d_{n-1} + (-1)^n} \to \frac{1}{n+1}$.  In the hats problem, this says that the chance the last hat is returned to its owner is approximately $\frac{1}{n+1}$ if all the previous hats have been returned incorrectly.  In fact, 
 it follows that if there have been exactly $d$ hats returned correctly among the first $(n-1)$ people, the last person receives the correct hat with probability approximately equal to $\frac{1}{n-d+1}$ when $d$ is far smaller than $n$.

\begin{figure}[h]
\tiny
\begin{tabular}{rccccccccc}
$n=1$:&    &    &    &    &  $f(1,0,0)$\\\noalign{\smallskip\smallskip}
$n=2$:&    &    &    &  $f(2,0,0,)$ &    &  $f(2,1,0)$\\\noalign{\smallskip\smallskip}
$n=3$:&    &    &  $f(3,0,0)$ &    &  $f(3,1,0)$ &    &  $f(3,2,0)$\\\noalign{\smallskip\smallskip}
$n=4$:&    &  $f(4,0,0)$ &    &  $f(4,1,0)$ &    &  $f(4,2,0)$ &    &  $f(4,3,0)$\\\noalign{\smallskip\smallskip}
$n=5$:&  $f(5,0,0)$ &    &  $f(5,1,0)$ &    &  $f(5,2,0)$ &    &  $f(5,3,0)$ &    &  $f(5,4,0)$\\\noalign{\smallskip\smallskip}
\end{tabular}
\caption{Equivalent triangle to Figure~\ref{fign=5}, but with $d=0$ by recurrence.   \label{fign=6} }

\end{figure}
\normalsize


\section{Partial Derangements}
\label{dZero}

In this section we let $d=0$, and focus on permutations with no fixed points in the first $k$ points.  
We proceed by proving a sequence of lemmas, then give an approximation for $f(n,k,0)$ and show that $p(n,k,0)$ and $f(n,k,0) $ are monotone in each of its variables, with a few noted exceptions. 

The functions will be then considered with some fixed points in Section~\ref{secResults} using Proposition~\ref{proprecursion}.
To start we will show that, for $n>3$, $f$ is decreasing in $k$ whenever it is also decreasing in $n$.

\begin{lemma}\label{lemX}
Let $n > 3$ and $0 < k < n$.  Then $f(n,k-1,0) <  f(n-1,k-1,0)$ if and only if $f(n,k,0) < f(n,k-1,0)$.
\end{lemma}

\begin{proof}

First let $k=1$.  We need to show that  $f(n,0 ,0) <  f(n-1,0 ,0)$ if and only if $f(n,1,0) < f(n,0 ,0)$.
In Table~\ref{tab3}, we have seen that $f(n,0,0) = \frac{1}{n}$, and thus $f(n,0 ,0) <  f(n-1,0 ,0)$.  We have also seen that $f(n,1,0) = \frac{n-2}{(n-1)^2}$.  Therefore, it is sufficient to show that $ \frac{n-2}{(n-1)^2} < \frac{1}{n}$, which is true for all $n$.

Now let $n > 3$ and $1 < k < n$, and consider the following events:
$$ U = \{ \alpha \in S_n \ | \ \alpha(1) \neq 1 \}$$
$$ V = \{ \alpha \in S_n \ | \ \alpha(i) \neq i \textrm{ for all } i \in [2,\ldots, k] \}$$
$$ W = \{ \alpha \in S_n \ | \ \alpha(k+1) \neq k+1 \}.$$

For any events $U$, $V$, and $W$ the following holds:
$$P(U | V) = P(U \cap W | V) + P( U \cap W^c | V) =
P(W |V)P(U | W \cap V) + P(W^c | V)P(U | W^c \cap V).$$
Thus $P(U | V)$ is always a convex combination (weighted average) of $P(U | W \cap V)$ and $P(U | W^c \cap V).$ Therefore either
$$ P(U | W \cap V) \leq  P(U | V) \leq P(U | W^c \cap V)  \textrm{ or }$$
$$ P(U | W \cap V) >  P(U | V) > P(U | W^c \cap V)  .$$

For our particular $U$, $V$, and $W$, $P(U | V) = 1 - f(n,k-1,0)$ and $P(U | V \cap W) = 1 - f(n,k,0)$. 
So either
$$  1 - f(n,k,0) \leq  1 - f(n,k-1,0) \leq P(U | W^c \cap V)  \textrm{ or }$$
$$  1 - f(n,k,0) >  1 - f(n,k-1,0)  > P(U | W^c \cap V)  .$$
Furthermore, 

\begin{align*}
P(U | W^c \cap V)  &  =  \frac{P(U \cap V \cap W^c)}{P(V\cap W^c)}  \\
& = \frac{P(W^c)P(U \cap V | W^c)}{P(W^c)P(V  | W^c)} \\
& =  \frac{P(U \cap V | W^c)}{P(V  | W^c)}   \\
& =   \frac{p(n-1,k,0)}{p(n-1,k-1,0)}       \textrm{ (by Corollary~\ref{cor_kij}) } \\
& = 1 - f(n-1,k-1,0)   \textrm{ (by Equation~\ref{eqftop}) }
\end{align*}

Thus either
\begin{align*}
1 - f(n,k,0) \leq  1 - f(n,k-1,0) \leq 1 - f(n-1,k-1,0)   \textrm{ or } \\
  1 - f(n,k,0) >  1 - f(n,k-1,0)  > 1 - f(n-1,k-1,0).
\end{align*}

\end{proof}

\begin{proposition} \label{prop1}
For all $0 < k \leq n$, $$p(n,k,0) \leq p(n,k-1,0).$$
\end{proposition}
\begin{proof}
This follows from the fact that $S_{n,k,0}$ a subset of the $S_{n,k-1,0}$.
\end{proof}

We point out that this claim is false for $p(n,k,d)$ when $d \neq 0$.

\begin{proposition} \label{prop3}
For $k < n$, $$p(n-1,k,0) \leq p(n,k,0).$$
\end{proposition}
\begin{proof}

By Proposition \ref{prop1} and Lemma~\ref{lem2}, 
\begin{align*}
p(n,k,0) =  \frac{k}{n} p(n-1,k-1,0) + (1-\frac{k}{n})p(n-1,k,0) \\
 \geq \frac{k}{n} p(n-1,k,0) + (1-\frac{k}{n}) p(n-1,k,0) = p(n-1,k,0).
 \end{align*}
\end{proof}

\begin{lemma}\label{lem4}
For $0 < k < n$, $$p(n,k,0) \leq p(n-1,k-1,0).$$
\end{lemma}
\begin{proof}

Again, by Proposition \ref{prop1} and Lemma \ref{lem2}, 

\begin{align*}
 p(n,k,0) &=  \frac{k}{n} p(n-1,k-1,0) + (1-\frac{k}{n})p(n-1,k,0)  \\
& \leq \frac{k}{n} p(n-1,k-1,0) + (1-\frac{k}{n}) p(n-1,k-1,0) = p(n-1,k-1,0).
 \end{align*}

\end{proof}
 
Note that this lemma does not hold for $k=n$ when $n$ is even, since

 $$p(n,n,0) - p(n-1,n-1,0) = \frac{(-1)^n}{n!}.$$

\begin{lemma}\label{lem5}
For $k < n$, $$p(n,k+1,0) = p(n,k,0) - \frac{1}{n} p(n-1,k,0).$$

\end{lemma}
\begin{proof}

To find this equality, we break $p(n,k,0)$ into two cases, based on whether $k+1$ is fixed.
If $k+1$ is not fixed, then there are exactly 0 fixed points in the first $k+1$ points.  This happens with probability $p(n,k+1,0)$.
If $k+1$ is fixed, then there 0 fixed points in the first $k$ points, and we know the image of $k+1$.  By Corollary~\ref{cor_kij}, this happens with probability $c(n-1,k,0)/n!$.   Thus
\begin{align*}
p(n,k,0): &= p(n,k+1,0) + \frac{ c(n-1,k,0)}{n!} \\
&=p(n,k+1,0)+ \frac{1}{n} p(n-1,k,0).
\end{align*}
\end{proof}

\begin{lemma}\label{lem6}
For $0 < k < n $, and $n \geq 3$, 
$$\frac{p(n,k+1,0)}{p(n,k,0)} = 1-\frac{1}{n} + \frac{k}{n^2(n-1)} \frac{p(n-2,k-1,0)}{p(n,k,0)}.$$

\end{lemma}
\begin{proof}

By Lemma \ref{lem5} and Lemma \ref{lem2}, 
\begin{align*}
p(n,k+1,0) & = p(n,k,0) - \frac{1}{n} p(n-1,k,0) \\
 & = p(n,k,0)(1-\frac{1}{n}) + \frac{1}{n}p(n,k,0) - \frac{1}{n} p(n-1,k,0)  \\
 & = p(n,k,0)(1-\frac{1}{n}) + \frac{1}{n} [p(n,k,0) - p(n-1,k,0)] \\
 & = p(n,k,0)(1-\frac{1}{n}) + \frac{1}{n} [ \frac{k}{n} p(n-1,k-1,0) + (1-\frac{k}{n})p(n-1,k,0) - p(n-1,k,0)]\\
 & = p(n,k,0)(1-\frac{1}{n}) + \frac{1}{n} \frac{k}{n} [p(n-1,k-1,0) - p(n-1,k,0)].\\
 \end{align*}
Then by Lemma \ref{lem5}, with $n$ and $k$ replaced by $n-1$ and $k-1$, we have $p(n-1,k,0)=p(n-1,k-1,0) - \frac{1}{n-1}p(n-2,k-1,0)$, and thus
$$ p(n,k+1,0) = p(n,k,0)[1-\frac{1}{n} + \frac{1}{n}\frac{k}{n} \frac{1}{(n-1)} p(n-2,k-1,0)].$$
Dividing through by $p(n,k,0)$ completes the proof.
\end{proof}

\begin{corollary}\label{cor1}
For $0 < k < n$,  $$\frac{p(n,k+1,0)}{p(n,k,0)} \geq  \frac{n-1}{n}.$$ Furthermore,
$$f(n,k,0) \leq \frac{1}{n}.$$

\end{corollary}

\begin{theorem}\label{thm1}
For $n>3$ and $k < n$, $f(n,k,0) \approx \frac{1}{n} - \frac{k}{n^2(n-1)}$, in particular:
$$ \frac{1}{n} - \frac{k}{n^2(n-1)}\frac{n}{(n-1)} \leq f(n,k,0) \leq
        \frac{1}{n} - \frac{k}{n^2(n-1)} \frac{(n-3)}{(n-2)}  .$$
\end{theorem}
\begin{proof}

In Table~\ref{tab3} we saw that $f(n,0,0) =  \frac{1}{n}$ and $f(n,1,0) =  \frac{(n-2)}{(n-1)^2}$ and thus the theorem is true for $k=0,1$.  Now assume $2 \leq k < n$.

By Lemma \ref{lem6}, Proposition \ref{prop3}, and Corollary \ref{cor1},
\begin{align*}
1-f(n,k,0) = \frac{p(n,k+1,0)}{p(n,k,0)} & = 1-\frac{1}{n} + \frac{k}{n^2(n-1)} \frac{p(n-2,k-1,0)}{p(n,k,0)}\\
& \leq 1-\frac{1}{n} + \frac{k}{n^2(n-1)} \frac{p(n,k-1,0)}{p(n,k,0)} \\
&  \leq 1-\frac{1}{n} + \frac{k}{n^2(n-1)} \frac{n}{(n-1)}
\end{align*}

By Lemma \ref{lem6}, Lemma \ref{lem4}, and Corollary \ref{cor1},
\begin{align*}
1-f(n,k,0) = \frac{p(n,k+1,0)}{p(n,k,0)} & = 1-\frac{1}{n} + \frac{k}{n^2(n-1)} \frac{p(n-2,k-1,0)}{p(n,k,0)}\\
& \geq 1-\frac{1}{n} + \frac{k}{n^2(n-1)} \frac{p(n-1,k-1,0))}{p(n-2,k-2,0)} \\
& \geq 1-\frac{1}{n} + \frac{k}{n^2(n-1)} \frac{n-3}{n-2}.
\end{align*}
\end{proof}

\begin{theorem}\label{thm2}
For $0 < k < n$ and $n \neq 3$, $f(n,k,0) < f(n,k-1,0)$.
\end{theorem}
\begin{proof}
By Lemma~\ref{lemX}, we can instead prove that $f(n,k-1,0) < f(n-1,k-1,0)$.

To start, Table~\ref{tab2} shows that this is true for all $n \leq 6$, and thus we assume that $n > 6$.

By Theorem \ref{thm1}, 
\begin{align*}
f(n,k-1,0) \leq \frac{1}{n} - \frac{(k-1)(n-3)}{n^2(n-1)(n-2)},\textrm{ and also } \\
\frac{1}{n-1} - \frac{(k-1)(n-1)}{(n-1)^2(n-2)^2} \leq f(n-1,k-1,0)
\end{align*}

and thus it suffices to show that 
\begin{align*}
\frac{1}{n-1} - \frac{1}{n}  >  \frac{(k-1)}{(n-1)(n-2)^2} - \frac{(k-1)(n-3)}{n^2(n-1)(n-2)}. \\
\end{align*}

This can be simplified as follows.
\begin{align*}
 \frac{1}{n-1}  - \frac{1}{n}  >  \frac{(k-1)}{(n-1)(n-2)^2}  - \frac{(k-1)(n-3)}{n^2(n-1)(n-2)}  \\
\frac{1}{n(n-1)} > \frac{(k-1)}{(n-1)(n-2)} \left[ \frac{1}{(n-2)} - \frac{(n-3)}{n^2} \right] \\
\frac{1}{n(n-1)} > \frac{(k-1)}{(n-1)(n-2)} \frac{5n-6}{n^2(n-2)} \\
1 > \frac{(k-1)(5n-6)}{n(n-2)^2}.
\end{align*}

By definition, $k \leq n-1$, so it is sufficient to show
$$1 > \frac{(n-2)(5n-6)}{n(n-2)^2}
 $$
which holds for all $n>6$.  
\end{proof}


\section{Proof of General Theorems}
\label{secResults}

In this section we now give the general results about 
$f(n,k,d)$, allowing $d$ to no longer be zero.  We utilize the results of Section~\ref{dZero}, and the recursion in Proposition~\ref{proprecursion}.

\begin{theorem}\label{thmnew1}
For $k < n$, $f(n,k,d) \approx \frac{1}{(n-d)} - \frac{(k-d)}{(n-d)^2(n-d-1)}$, in particular:
$$\frac{1}{(n-d)} - \frac{(k-d)}{(n-d)^2(n-d-1)} \frac{(n-d)}{(n-d-1)} \leq  f(n,k,d) \leq               \frac{1}{(n-d)} - \frac{(k-d)}{(n-d)^2(n-d-1)} \frac{(n-d-3)}{(n-d-2)}    .$$\end{theorem}
\begin{proof}
This follows from Theorem~\ref{thm1}, and Proposition~\ref{proprecursion}.
\end{proof}

\begin{theorem}\label{thmnewk}
Let $0 \leq i < j < n$.
If $(n-d) \neq 3$, then $f(n,j,d) < f(n,i,d)$.   Furthermore, $f(n,j,d) > f(n,i,d)$ if and only if $n=j+1=i+2=d+3$.    
\end{theorem}
\begin{proof}
By Theorem~\ref{thm2} and Proposition~\ref{proprecursion}, $f(n,j,d) < f(n,i,d)$ for $0 < i < j < n$ as long as $(n-d) \neq 3$.  
When $(n-d)=3$, by Proposition~\ref{proprecursion} we can consider $f(3,k,0)$.  Here $f(3,2,0) = \frac{1}{3}$, $f(3,1,0) = \frac{1}{4}$, and $f(3,0,0) = \frac{1}{3}$.    
Thus $f(3,j,0) \leq f(3,i,0)$ for $0 < i < j < n$, except when $j=2$ and $i=1$.   Using Proposition~\ref{proprecursion} a final time, we have our result.
\end{proof}

\begin{theorem}\label{thmnewn}
For $0 \leq k < n$ and $m >n$, $f(m,k,d) < f(n,k,d)$ except when $n=d+2$ and $k=d+1$.
In particular when $d=0$, $f(n,k,0)$ is decreasing in $n$ when $k \neq 1$.
\end{theorem}
\begin{proof}
The exception to the decreasing behavior can be seen in Table~\ref{tab2} when $f(2,1,0) < f(m,1,0)$ for $m>2$.  The general case follows from Theorem~\ref{thmnewk}, Proposition~\ref{proprecursion}, and Lemma~\ref{lemX}.
\end{proof}

\begin{theorem} \label{thmnewd}
Let $0 \leq c < d \leq k < n$.  Then $f(n,k,d) > f(n,k,c)$ except when $n=k+1=d+2$; in this case $f(n,n-1,n-2)=0$.
\end{theorem}

\begin{proof}
The statement can be checked when $n-d < n-c \leq 6$ using Proposition~\ref{proprecursion} and Table~\ref{tab2}.
Now assume $n-c > 6$.  We will show that $f(n,k,c+1) > f(n,k,c)$ for all $c$.  Using Proposition~\ref{proprecursion}, this is equivalent to showing $f(n-c-1,k-c-1,0) > f(n-c,k-c,0)$.
The rest of the proof then follow using the inequalities of Theorem~\ref{thm1} in the same way as in the proof of Theorem~\ref{thm2}; the full details are omitted.
\end{proof}


\section{Acknowledgement}
We would like to thank Peter Doyle and Larry Shepp for helpful comments on an early draft of this paper.

\bibliographystyle{acm}
 \bibliography{FixedPointsBib}
 
\end{document}